\documentclass[12pt]{article}
\usepackage{latexsym,amsmath,amssymb}
\begin{document}

\noindent J Theor Probab (2015) 28:1337-1353\vspace{0.8cm}\\
\noindent
{\bf \large Moment Determinacy of Powers and Products\\
of Nonnegative Random Variables}

\vspace{0.8cm} \noindent {\bf Gwo Dong Lin$^{1}$ \ $\bullet$ \
Jordan Stoyanov$^{2}$ }

\vspace{0.2cm}\noindent $^{1}$ Institute of Statistical Science,
Academia Sinica, Taipei 11529, Taiwan (ROC) \\e-mail:
gdlin@stat.sinica.edu.tw

\vspace{0.2cm}\noindent $^{2}$ School of Mathematics $\&$
Statistics, Newcastle University, Newcastle upon Tyne NE1 7RU,
UK\\
e-mail: stoyanovj@gmail.com

\vspace{1.0cm}\noindent {\bf Abstract:} We find conditions which
guarantee moment (in)determinacy of powers and products of
nonnegative random variables. We establish new and general results
which are based either on the rate of growth of the moments  of a
random variable or on  conditions about the distribution itself.
For the class of generalized gamma random variables we show that
the power and the product of such variables
share the same moment determinacy property. A similar statement
holds for half-logistic random variables. Besides answering new
questions in this area, we  either extend some previously known
results or provide new and transparent proofs of  existing
results.

\vspace{0.2cm}\noindent

\vspace{0.2cm}\noindent

\vspace{0.4cm}\noindent {\bf 1. Introduction} \

\vspace{0.1cm}\noindent Throughout the paper we assume that $\xi$
is a nonnegative random variable defined on a given probability
space $(\Omega, {\cal F}, {\bf P})$ with finite moments ${\bf
E}[\xi^k], \ k=1,2,\ldots.$ Let further $\xi_1,\xi_2, \ldots,
\xi_n$ be independent copies of $\xi.$ Here $n\ge 1$ is a fixed
integer number. Of interest to us are the following two random
variables, the power and the product:
\[
X_n=\xi^n \ \mbox{ and } \ Y_n=\xi_1\xi_2 \cdots \xi_n.
\]

Each of the variables $X_n$ and $Y_n$ also has all moments finite. Thus
the natural problem arising here is to study, characterize and compare the moment (in)determinacy
of these two variables.
Since $X_n$ and $Y_n$ take values in ${\mathbb R}^+=
[0,\infty)$, this means that we deal with the Stieltjes moment problem.

We find conditions on $\xi$ and $n$ guaranteeing that
$X_n$ and $Y_n$ are M-determinate (uniquely
determined by the moments), and other conditions when they are
M-indeterminate (nonunique in terms of the moments).
In these two cases we use the abbreviations M-det and M-indet
for both random variables and their distributions.

In our reasonings we use classical or new conditions such as
Cram\'{e}r's condition, Carleman's condition, Hardy's condition and Krein's condition.
 The reader may find it useful to consult available sources,
 among them are \cite{l1997}, \cite{ls2009},  \cite{p2001}, \cite{p1998},
 \cite{s2000}, \cite{s2013}
 and \cite{sl2012}.
For reader's convenience, we have included these conditions when formulating
our results.

To study powers, products, etc., or other nonlinear
transformations of random data (called sometimes Box-Cox
transformations), is a challenging probabilistic problem which is
of independent interest. Note however that products and powers of
random variables considered in this paper and the results
established are definitely related to contemporary stochastic
models of real and complex phenomena; see, e.g., \cite{cm1994},
\cite{da2010}, \cite{gs2004} and \cite{pb2010}.

In this paper we deal with new problems and present new results
with their proofs. We establish new and general criteria which are
then applied to describe the moment (in)determinacy of the above
 transformations. We also provide new proofs of
some known results with reference to the original papers. Our
results complement previous studies or represent different aspects
of existing studies on this topic; see, e.g., \cite{b2005},
\cite{d1997}, \cite{lh1997},
  \cite{os2010},   \cite{pb2010} and \cite{s2000}.

The approach and the results in this paper can be further extended
to distributions on the whole real line (Hamburger moment problem,
see \cite{sld2014}). Also, they can be used to
characterize the moment determinacy properties of nonlinear
transformations of some important sub-classes of distributions
such as, e.g., the subexponential distributions; see \cite{fkz2011}.

The material is divided into relatively short sections each
dealing with a specific question related to a general or specific distribution.
General results are included in Sections 2, 4, 6, 7 and 9. Sections 3, 5,  8 and 10
deal with powers and products based on the generalized gamma distribution while Section
11 is based on half-logistic distribution.
All statements are followed by
detailed proofs.

\vspace{0.5cm}\noindent {\bf 2. Comparing the moment determinacy of powers and
products}

\vspace{0.1cm}
The power $X_n=\xi^n$ and  the product $Y_n=\xi_1\xi_2 \cdots \xi_n$ have some
`similarity'. They both are defined in terms of $n$ and $\xi$ or of $n$ independent copies of $\xi$,
and both have all moments finite. Thus we arrive naturally to the question:

 Is it true that the random variables $X_n$ and $Y_n$ share the same moment determinacy property?

 If the generic random variable $\xi$ has a bounded support, then so does each of
$X_n$ and $Y_n$, and hence both $X_n$ and $Y_n$ have all moments
finite and both are M-det. This simple observation shows that interesting is to
study powers and products based on a
random variable $\xi$ with unbounded support contained in
${\mathbb R}^+$ and such that $\xi$ has all moments finite.

Let us mention first a special case. Suppose $\xi \sim Exp(1),$  the
standard exponential distribution. Then the power $X_n=\xi^n$ is
M-det iff the product $Y_n=\xi_1\xi_2\cdots\xi_n$ is M-det and this is true iff
$n\leq 2$ (see, e.g., \cite{b2005} and 
\cite{os2010}). This means that for any $n=1,2,\ldots,$ the power
$X_n$ and the product $Y_n$ share the same moment determinacy
property. Since Weibull random variable is just a power of the
exponential one, it follows that if $\xi$ obeys a Weibull
distribution, then for any $n=1,2,\ldots,$ the power $X_n$ and the
product $Y_n$ also have the same moment determinacy property.
Therefore, the answer to the above question is positive for at
least some special distributions including Weibull distributions.
In this paper we will explore more distributions (see Theorem 6
and Section 11 below).

Note that in general,  we have, by Lyapunov's inequality,
\begin{eqnarray}
{\bf E}[X_n^s]={\bf E}[\xi^{ns}] \geq ({\bf E}[\xi^s])^n={\bf
E}[Y_n^s]~~\hbox{for all real}~ s>0.
\end{eqnarray}
We use this moment inequality to establish a result which involves
three of the most famous conditions  for moment determinacy
(Carleman's, Cram\'{e}r's and Hardy's). For more details about
Hardy's condition, see 
\cite{sl2012}.

\vspace{0.2cm} \noindent {\bf Proposition 1} \ (i) {\it If $X_n$
satisfies Carleman's condition (and hence is M-det), i.e.,
$\sum_{k=1}^{\infty}({\bf E}[X_n^k])^{-1/(2k)}=\infty,$ then so does $Y_n$.} \\
(ii) {\it If $X_n$ satisfies Cram\'er's condition (and hence is
M-det),
i.e., ${\bf E}[\exp({cX_n})]<\infty$ for some constant $c>0,$ then so does $Y_n$.}\\
(iii) {\it If $X_n$ satisfies Hardy's condition (and hence is
M-det), i.e., ${\bf E}[\exp({c\sqrt{X_n}})]<\infty$ for some
constant $c>0,$ then so does $Y_n$.}

\vspace{0.1cm}\noindent
{\it Proof} \ Part (i) follows immediately from (1). Parts (ii) and
(iii) follow from the fact that for each real $s>0$,
$${\bf E}[\exp(cX_n^s)]=\sum_{k=0}^{\infty}\frac{c^k}{k!}{\bf E}
[(X_n^s)^{k}]\geq \sum_{k=0}^{\infty} \frac{c^k}{k!}{\bf
E}[(Y_n^s)^{k}]={\bf E}[\exp(cY_n^s)].$$
\hspace{\fill}$\Box$

\vspace{0.2cm}\noindent {\bf Corollary 1} \ {\it If $\xi$
satisfies Cram\'er's condition and if $n=2$, then both $X_2=\xi^2$
and $Y_2=\xi_1\xi_2$ are M-det, and hence $X_2$ and $Y_2$ share
the same moment determinacy property.}

\vspace{0.1cm}\noindent {\it Proof} \ Note that $\xi$ satisfies
Cram\'er's condition iff $X_2$ satisfies Hardy's condition. Then
by Proposition 1(iii), both $X_2$ and $Y_2$ are M-det as claimed
above. \hspace{\fill}$\Box$

\vspace{0.5cm}\noindent {\bf 3. Generalized gamma distributions.
Part (a)}

\vspace{0.1cm} Some of our results can be well illustrated if we
assume that the generic random variable $\xi$ has a generalized
gamma distribution. We write $\xi \sim GG(\alpha, \beta, \gamma)$
if $\xi$ has the following density function $f$:
$$f(x)=cx^{\gamma-1}e^{-\alpha
x^{\beta}},~x\geq 0,$$ where $\alpha, \beta, \gamma>0$, $f(0)=0$
if $\gamma\ne 1$, and
$c=\beta\alpha^{\gamma/\beta}/\Gamma(\gamma/\beta)$ is the norming
constant. Note that $GG(\alpha, \beta, \gamma)$ is a rich class
containing several commonly used distributions such as
exponential, Weibull, half-normal and chi-square.

It is known that the power $X_n=\xi^n$ is M-det iff
 $n\leq 2\beta$ (see, e.g.,  \cite{pk1992} and \cite{s2013}).
We claim now  that for $n\leq 2\beta$, the product $Y_n=\xi_1\xi_2\cdots\xi_n$
is also
M-det. To see this, we note first that the density function $h_n$
of the random variable
$\sqrt{X_n}$ is
$$h_n(z)=\frac{2c}{n}z^{2\gamma/n-1}e^{-\alpha
z^{2\beta/n}},~z\geq 0.$$ This in turn implies that $X_n$
satisfies Hardy's condition if $2\beta/n\geq 1$, hence so does
$Y_n$ for $n\leq 2\beta$ by Proposition 1(iii).

To obtain further results, it is quite useful to write the explicit
form of the density of the product $Y_2=\xi_1\xi_2$ when
$\xi$ has the  generalized gamma distribution. This involves the function $K_0(x), \
x >0$, the modified Bessel function of the second kind. Its
definition and approximation  are given as follows:
\begin{eqnarray*}
K_0(x)&=& \frac{1}{2}\int_0^{\infty}t^{-1}e^{-t-x^2/(4t)}dt,~~x>0, \\
&=&\left(\frac{\pi}{2x}\right)^{1/2}e^{-x}\left[1-\frac{1}{8x}
\left(1-\frac{9}{16x}\left(1-\frac{25}{24x}\right)\right)+o(x^{-3})\right]
~~\hbox{as}~x\rightarrow\infty
\end{eqnarray*}
(see, e.g., 
\cite{gkk2012} and 
\cite{malh2005}, pp. 37--38).

\vspace{0.2cm}\noindent
{\bf Lemma 1} (See also 
\cite{mali1968})\ {\it Let $Y_2=\xi_1\xi_2$, where $\xi_1$ and
$\xi_2$ are independent random variables having the same
distribution $GG(\alpha, \beta, \gamma).$
Then the
density function $g_2$ of $Y_2$ is
\begin{eqnarray*}
g_2(x)&=&
\frac{2c^2}{\beta}x^{\gamma-1}K_0\left(2\alpha x^{\beta/2}\right),~~x>0,\\
&\approx& Cx^{\gamma-\beta/4-1}e^{-2\alpha
x^{\beta/2}},~~\hbox{as}~x\rightarrow\infty.
\end{eqnarray*}
} \vspace{0.1cm}\noindent {\it Proof} (Method I) Let $G_2$ be the
distribution function of $Y_2$. Then
$$\overline{G}_2(x):=1-G_2(x)={\bf P}[Y_2>x]=\int_0^{\infty}{\bf P}[\xi_1>x/y]
cy^{\gamma-1}e^{-\alpha y^{\beta}}dy,~~x>0,$$ and hence the
density of $Y_2$ is
\begin{eqnarray*}
g_2(x)&=&c^2x^{\gamma-1}\int_0^{\infty}y^{-1}e^{-\alpha
x^{\beta}/y^{\beta} -\alpha y^{\beta}}dy
=\frac{c^2}{\beta}x^{\gamma-1}\int_0^{\infty}t^{-1}e^{-t-(\alpha^2
x^{\beta})/t}dt\\
&=&\frac{2c^2}{\beta}x^{\gamma-1}K_0\left(2\alpha
x^{\beta/2}\right),~~x>0.
\end{eqnarray*}
(Method II) We can use the moment function (or Mellin transform)
usually denoted by $\cal M$,
because it uniquely determines the corresponding distribution. To
do this, we note that
$${\cal M}(s) =: {\bf E}[Y_2^s]=({\bf E}[\xi_1^s])^2, ~~{\bf
E}[\xi_1^s]=c\Gamma((\gamma+s)/\beta)\left(\beta\alpha^{(\gamma+s)/
\beta}\right)^{-1},~~\hbox{and}$$
$$\int_0^{\infty}x^sK_0(x)dx=2^{s-1}(\Gamma((s+1)/2))^2~~
\hbox{for all}~~s>0$$ (see, e.g., \cite{gr2007}, p. 676,
Formula 6.561(16)). We omit the detailed
calculation. \hspace{\fill} $\Box$

\vspace{0.2cm} It may look surprising, but it is well-known, that
several commonly used distributions are related to the Bessel
function in such a natural way as in Lemma 1. For example, if
$\xi$ is a half-normal random variable, i.e., $\xi \sim
GG(\frac{1}{2}, 2, 1)$ with the density
$f(x)=\sqrt{2/\pi}e^{-x^2/2},~x\geq0$, then $Y_2$ has the density
function $g_2(x)=(2/{\pi})K_0(x)\approx
C_2x^{-1/2}e^{-x}~~\hbox{as}~x\rightarrow\infty$, with the moment
function  ${\cal M}(s)={\bf
E}[Y_2^s]=(2^s/\pi)\Gamma^2((s+1)/2),~s>-1$. The distribution of
$Y_2=\xi_1\xi_2$ may be called the half-Bessel distribution and
its symmetric counterpart with density
$h_2(x)=(1/{\pi})K_0(x),~x\in {\mathbb R}= (-\infty,\infty)$, is
called the standard Bessel distribution. Note that $K_0$ is an
even function and $h_2$ happens to be the density of the product
of two independent standard normal random variables; see also
\cite{d1997}.
  It can be checked  that for real $s>0$ we have $({\bf
E}[(Y_2^{s})^n])^{-1/(2n)}\approx C_sn^{-s/2}$ as
$n\rightarrow\infty$, and hence   $Y_2^s$ satisfies Carleman's
condition iff $s\leq 2$. Actually, it follows from the density
 $g_2$ and its asymptotic behavior that $Y_2$ satisfies
Cram\'er's condition. Therefore, by Hardy's criterion, the square
of $Y_2$, i.e., $Y_2^2=\xi_1^2\xi_2^2$, is M-det.

Let us express the latter by words: The square of the product of
two independent half-normal random variables is M-det. Since $\xi^2 = \chi^2_1$, we
conclude also that the product of two independent
$\chi^2$-distributed random variables is M-det. In addition, these
properties can be compared with the known fact that the power 4 of
a normal random variable is M-det (see, e.g,
\cite{b1988} or \cite{s2000}).

\vspace{0.5cm} \noindent {\bf 4. Slow growth rate of the moments
implies moment determinacy}

\vspace{0.1cm} It is known and well understood that the moment
determinacy of a distribution depends on the rate of growth of the
moments. Let us establish first  results which are of a general
and independent interest. Later  we will apply them and make
conclusions about powers and products of random variables.

Suppose $X$ is a nonnegative random variable with finite moments
$m_k={\bf E}[X^k]$, $k=1,2,\ldots.$ To avoid trivial cases, we
assume that $m_1>0$, meaning that $X$ is not a degenerate random
variable at $0$.

\vspace{0.2cm}\noindent
{\bf Lemma 2} \ {\it For each $k\geq 1$, we have the following properties:}\\
(i) {\it $m_k\leq m_{k+1}$ if $m_1\ge 1$, and} \\
(ii) {\it $m_1m_k\leq m_{k+1}.$}

\vspace{0.1cm}\noindent {\it Proof} \ By Lyapunov's inequality, we
have $(m_k)^{1/k}\leq (m_{k+1})^{1/(k+1)}.$ Therefore,
$$\frac{1}{k}\log m_k\leq \frac{1}{k+1}\log m_{k+1}\leq
\frac{1}{k}\log m_{k+1}, ~\hbox{if} ~m_1\ge 1,$$ and hence
$m_k\leq
m_{k+1}$ if $m_1\ge 1$. This proves claim (i).\\
To prove claim (ii), we use the relations $m_1\leq (m_k)^{1/k}\leq
(m_{k+1})^{1/(k+1)},$ implying that
$$m_1m_k\leq (m_k)^{1/k}m_k=m_k^{(k+1)/k}\leq m_{k+1}.$$
\hspace{\fill}$\Box$

\vspace{0.2cm} In Lemma 2, claim (i) tells us that the moment
sequence $\{m_k, k=1,2,\ldots \}$ is nondecreasing if $m_1\ge 1$,
while claim (ii) shows that the ratio $m_{k+1}/m_k$ has a lower
bound $m_1$ whatever the nonnegative random variable $X$ is. The
next theorem provides the upper bound of the ratio $m_{k+1}/m_k$,
or, equivalently, of the growth rate of the moments $m_k$ for
which $X$ is M-det.

\vspace{0.2cm}\noindent {\bf Theorem 1} \ {\it Let
$m_{k+1}/m_k={\cal O}((k+1)^2)$ as $k\rightarrow\infty$. Then $X$
satisfies Carleman's condition and is M-det. (We refer to the
constant 2, the exponent in the term ${\cal O}((k+1)^2)$, as the
rate of growth of the moments of $X$.) }

\vspace{0.1cm}\noindent {\it Proof} \ By the assumption, there
exists a constant $C>0$ such that
$$m_k^{(k+1)/k}\leq m_{k+1}\leq C(k+1)^2m_k~~\hbox{for all
large}~k,$$ which implies
$$m_k^{1/k}\leq  C(k+1)^2~~\hbox{for all
large}~k,$$ and hence
$$m_k^{-1/(2k)}\geq  C^{-1/2}(k+1)^{-1}~~\hbox{for
large}~k.$$ Therefore, $X$ satisfies Carleman's condition
$\sum_{k=1}^{\infty}m_k^{-1/(2k)}=\infty$, and is M-det.
\hspace{\fill}$\Box$

\vspace{0.2cm} We can slightly extend Theorem 1 as follows. For a
real number $a$ we denote by $\lfloor a \rfloor$ the largest
integer which is less than or equal to $a$.

\vspace{0.2cm}\noindent {\bf Theorem 1$^{\prime}$} \ {\it Suppose
there is a real number $a \geq 1$ such that the moments of the
random variable $X$ satisfy the condition $m_{k+1}/m_k={\cal
O}((k+1)^{2/a})$ as $k\rightarrow\infty$. Then the power
$X^{\lfloor a \rfloor}$ satisfies Carleman's condition and is
M-det.}

\vspace{0.1cm}\noindent {\it Proof} \ Note that
\begin{eqnarray*}\frac{{\bf E}[(X^{\lfloor a
\rfloor})^{k+1}]}{{\bf E}[(X^{\lfloor a \rfloor})^k]}&=&\frac{{\bf
E}[X^{{\lfloor a \rfloor}k+{\lfloor a \rfloor}}]}{{\bf
E}[X^{{\lfloor a \rfloor}k+{\lfloor a \rfloor}-1}]}\frac{{\bf
E}[X^{{\lfloor a \rfloor}k+{\lfloor a \rfloor}-1}]}{{\bf
E}[X^{{\lfloor a \rfloor}k+{\lfloor a
\rfloor}-2}]}\cdots\frac{{\bf E}[X^{{\lfloor a \rfloor}k+1}]}{{\bf
E}[X^{{\lfloor a \rfloor}k}]}\\&=&{\cal O}((k+1)^{(2/a)\lfloor a
\rfloor})={\cal
O}((k+1)^2)~~\hbox{as}~~k\rightarrow\infty.\end{eqnarray*} Hence,
by Theorem 1,  $X^{\lfloor a \rfloor}$ satisfies Carleman's
condition and is M-det. \hspace{\fill} $\Box$

\vspace{0.2cm}\noindent {\bf Theorem 2} \ {\it Let $\xi,~\xi_i,
~i=1,2,\ldots,n$, be defined as before and $Y_n=\xi_1\cdots\xi_n.$
If $\xi$ and the index $n$ are such that
\[
{\bf E}[\xi^{k+1}]/{\bf E}[\xi^k]={\cal O}((k+1)^{2/n}) \ \mbox{ as } k\rightarrow \infty,
\]
then $Y_n$ satisfies Carleman's condition and is M-det.}

\vspace{0.1cm}\noindent {\it Proof} \ By the assumption, we have
\[
{\bf E}[Y_n^{k+1}]/{\bf E}[Y_n^k]=({\bf E}[\xi^{k+1}]/{\bf E}[\xi^k])^n={\cal O}((k+1)^2)
\ \mbox{ as }~ k\rightarrow \infty.
\]
This, according to Theorem 1, implies the validity of Carleman's
condition for $Y_n,$ hence $Y_n$ is M-det as stated above.
\hspace{\fill}$\Box$

\vspace{0.2cm}\noindent {\bf Theorem 2$^{\prime}$} \ {\it Let
$a\geq 1$. If
\[
{\bf E}[\xi^{k+1}]/{\bf E}[\xi^k]={\cal O}((k+1)^{2/a}) \ \mbox{ as }
k\rightarrow \infty,
\]
then $Y_{\lfloor a \rfloor}$ satisfies Carleman's condition and is
M-det.}

\vspace{0.1cm}\noindent {\it Proof} \ Note that
\begin{eqnarray*}
&~&{\bf E}[Y_{\lfloor a \rfloor}^{k+1}]/{\bf E}[Y_{\lfloor a
\rfloor}^k]=({\bf E}[\xi^{k+1}]/{\bf E}[\xi^k])^{\lfloor a
\rfloor}\\ &=&{\cal O}((k+1)^{(2/a){\lfloor a \rfloor}})={\cal
O}((k+1)^2)
 \ \mbox{ as }
k\rightarrow \infty. \end{eqnarray*} The conclusions follow from
Theorem 1. \hspace{\fill} $\Box$

\vspace{0.5cm}\noindent {\bf 5. Generalized gamma distributions. Part (b) }

\vspace{0.1cm} We now apply the general results, Theorems 1 and 2
in Section 4, to give an alternative proof of the moment
determinacy established in Section 3.

Let, as before, $\xi \sim GG(\alpha, \beta, \gamma).$
We claim that for
$n\leq 2\beta$, both $X_n=\xi^n$ and $Y_n=\xi_1\xi_2\cdots\xi_n$ are M-det.
To see this, we first calculate that
$$
\frac{{\bf E}[X_n^{k+1}]}{{\bf E}[X_n^{k}]}=\frac{{\bf E}
[\xi^{n(k+1)}]}{{\bf E}[\xi^{nk}]}=
\frac{\Gamma((\gamma+n(k+1))/\beta)}{\alpha^{n/\beta}
\Gamma((\gamma+nk)/\beta)}\approx
(n/\alpha\beta)^{n/\beta} (k+1)^{n/\beta}~~ \hbox{as}~
k\rightarrow \infty.
$$
For this relation we have used the approximation of the gamma function:
$$\Gamma(x)\approx \sqrt{2\pi}x^{x-1/2}e^{-x}~~\hbox{
as}~x\rightarrow \infty$$ (see, e.g., \cite{ww1927}, p. 253).
 Then by Theorem 1, $X_n$ is M-det if $n\leq 2\beta$, and by Theorem 2,
 $Y_n$ is M-det if $1/\beta\leq 2/n$, i.e., if $n\leq 2\beta$, because
${\bf E}[\xi^{k+1}]/{\bf E}[\xi^k]={\cal O}((k+1)^{1/\beta})$ as $k\rightarrow \infty$.

For example, if $\xi \sim Exp(1) = GG(1, 1, 1)$, then the product $Y_2=\xi_1\xi_2$ is
M-det. In fact, by Lemma 1, the density  $g_2$ of $Y_2$ is
$g_2(x)=2K_0(2\sqrt{x})\approx Cx^{-1/4}e^{-2\sqrt{x}}$ as
$x\rightarrow\infty$, where $K_0$ is the modified Bessel function
of the second kind (see also \cite{mt1986}, p. 417, and \cite{gr2007}, p.
917, Formula 8.432(8)). If $\xi \sim GG(1/2, 2, 1)$, the
half-normal distribution, then $Y_n=\xi_1\xi_2\cdots\xi_n$ is
M-det for $n\leq 4$.  As mentioned before, the density function of
the product of two half-normals is $g_2(x)=(2/{\pi})K_0(x)\approx
C_2x^{-1/2}e^{-x}\hbox{ as }~x\rightarrow\infty$.

\vspace{0.5cm} \noindent {\bf 6. More results related to Theorems
1 and 2}

\vspace{0.1cm} Under the same assumption as that in Theorem 1, we
even have a stronger statement; see Theorem 3 below. Note that its
proof does not
 use Lyapunov's inequality, and that Hardy's condition
implies Carleman's condition. For convenience, we recall in the
next lemma a characterization of Hardy's condition in terms of the
moments (see \cite{sl2012}, Theorem 3).

\vspace{0.2cm} \noindent {\bf Lemma 3} \ {\it Let $a\in (0,1]$ and
let $X$ be a nonnegative random variable. Then ${\bf
E}[\exp({c{X}^a})]<\infty$ for some constant $c>0$ iff ${\bf
E}[X^k]\leq c_0^k\,\Gamma(k/a+1), ~k=1,2,\ldots,$ for some
constant $c_0>0$ (independent of $k$). In particular, $X$
satisfies Hardy's condition, i.e., ${\bf
E}[\exp({c\sqrt{X}})]<\infty$ for some constant $c>0$, iff ${\bf
E}[X^k]\leq c_0^k\,(2k)!, ~k=1,2,\ldots,$ for some constant
$c_0>0$ (independent of $k$).}

\vspace{0.2cm} \noindent {\bf Theorem 3} \ {\it Suppose $X$ is a
nonnegative random variable with finite moments $m_k={\bf E}[X^k],
\ k=1,2,\ldots$, such that the condition in Theorem 1 holds: \
$m_{k+1}/m_k = {\cal O}((k+1)^2)$ as $k \to \infty.$ Then $X$
satisfies Hardy's condition, and is M-det.}

\vspace{0.1cm} \noindent {\it Proof} \ By the assumption, there
exists a constant $c_*\geq m_1>0$ such that
\begin{eqnarray*}
m_{k+1}\leq c_*(k+1)^2m_k~~\hbox{for}~k=0,1,2,\ldots,
\end{eqnarray*}
where $m_0\equiv 1$. This implies that
\begin{eqnarray*}
m_{k+1}\leq (c_*/2)(2k+2)(2k+1)m_k~~\hbox{for}~k=0,1,2,\ldots,
\end{eqnarray*}
and hence  $m_{k+1}\leq
(c_*/2)^{k+1}\Gamma(2k+3)m_0~~\hbox{for}~k=0,1,2,\ldots.$ Taking
$c_0=c_*/2$,
$$m_{k+1}\leq c_0^{k+1}\Gamma(2k+3)~~\hbox{for}~k=0,1,2,\ldots,$$
or, equivalently,
$$m_{k}\leq c_0^{k}\Gamma(2k+1)~~\hbox{for}~k=1,2,\ldots.$$
Hence $X$ satisfies Hardy's condition by Lemma 3.  \hspace{\fill}
$\Box$

\vspace{0.2cm}\noindent {\bf Remark 1} \ The constant 2 (the
growth rate of the moments) in the condition of Theorem 1 is the
best possible in the following sense. For each $\varepsilon>0$,
there exists a random variable $X$ such that $m_{k+1}/m_k={\cal
O}( (k+1)^{2+\varepsilon})$ as $k\rightarrow\infty$, and $X$ is
M-indet. To see this, let us consider $X=\xi\sim GG(1, \beta, 1),$
which has density $f(x)=c\exp(-x^{\beta}), ~x>0.$ We have
$$
\frac{{\bf E}[\xi^{k+1}]}{{\bf E}[\xi^{k}]}=
\frac{\Gamma((k+2)/\beta)}{{\Gamma((k+1)/\beta)}}\approx
\beta^{-1/\beta}(k+1)^{1/\beta} ~~ \hbox{as}~ k\rightarrow
\infty.
$$
If for $\varepsilon>0$ we take $\beta
=\frac{1}{2+\varepsilon}<\frac12$, then ${{\bf E}[\xi^{k+1}]}/{{\bf
E}[\xi^{k}]}=$ ${\cal O}( (k+1)^{2+\varepsilon}) ~~ \hbox{as}~
k\rightarrow \infty.$ However, as mentioned before,  $X$ is M-indet.

\vspace{0.2cm}\noindent {\bf Remark 2} \ The constant $2/n$ in the
condition of Theorem 2 is the best possible. Indeed, we can show
that for each $\varepsilon>0$, there exists a random variable
$\xi$ such that ${\bf E}[\xi^{k+1}]/{\bf E}[\xi^k]={\cal
O}((k+1)^{2/n+\varepsilon})$ as $k\rightarrow \infty$, but $Y_n$
is M-indet. To see this, let us consider $X=\xi\sim GG(1, \beta,
1)$. For each $\varepsilon>0$, take $\beta=1/(2/n+\varepsilon)$,
then
$$
\frac{{\bf E}[\xi^{k+1}]}{{\bf E}[\xi^{k}]}=
\frac{\Gamma((k+2)/\beta)}{{\Gamma((k+1)/\beta)}}=
{\cal O}\left((k+1)^{2/n+\varepsilon}
\right)~~ \hbox{as}~ k\rightarrow \infty.
$$
However, since $n>2\beta$,  $Y_n$ is M-indet
(compare this with the statement in Section 10).

\vspace{0.5cm} \noindent {\bf 7. Faster growth rate of the moments
implies moment indeterminacy}

\vspace{0.1cm} We now establish a result which is converse to
Theorem 1.

\vspace{0.2cm}\noindent {\bf Theorem 4} \ {\it Suppose $X$ is a
nonnegative random variables whose moments $m_k, \ k=1,2,\ldots$,
are  such that $m_{k+1}/m_k\geq C(k+1)^{2+\varepsilon}$ for all
large $k$, where $C$ and $\varepsilon$ are positive constants.
Assume further that $X$ has a density  $f$ satisfying the
condition: for some $x_0>0$, $f$  is positive and differentiable
on $[x_0,\infty)$ and
\begin{eqnarray}L_{f}(x):=-\frac{xf^{\prime}(x)}{f(x)}\nearrow
\infty~~\hbox{as}~~x_0<x\rightarrow \infty.\end{eqnarray} Then $X$
is M-indet.}

\vspace{0.1cm}\noindent
{\it Proof} \ Without loss of generality we can
assume that $m_{k+1}/m_k\geq C(k+1)^{2+\varepsilon}$ for each $k\geq 1$. Therefore,
$$m_{k+1}\geq C^k ((k+1)!)^{2+\varepsilon}m_1~~\hbox{for}~k=1,2,\ldots.$$
Taking $C_0=\min\{C,m_1\}$, we have
$$m_{k+1}\geq C_0^{k+1} ((k+1)!)^{2+\varepsilon}~~\hbox{for}~k=1,2,\ldots,$$
or, equivalently,
$$m_{k}\geq C_0^{k} (k!)^{2+\varepsilon}=C_0^{k}
(\Gamma(k+1))^{2+\varepsilon}~~\hbox{for}~k=2,3,\ldots.
$$
Since $\Gamma(x+1)=x\Gamma(x)\approx
\sqrt{2\pi}\,x^{x+1/2}\,e^{-x}$ as $x\rightarrow \infty$,
we have that for some constant $c>0$,
$$
m_{k}^{-1/(2k)}\leq  C_0^{-1/2}
(\Gamma(k+1))^{-(2+\varepsilon)/(2k)}\approx
ck^{-1-\varepsilon/2}~~\hbox{for all large}~k.
$$
This implies that the Carleman quantity for the moments of $f$ is
finite:
$$
{\bf C}[f]:=
\sum_{k=1}^{\infty}m_k^{-1/(2k)}<\infty.
$$
We sketch the rest of the proof. Following the proof of Theorem 3
in \cite{l1997}, we first construct a symmetric distribution $G$
on ${\mathbb R}$, obeyed by a random variable $Y$, such that ${\bf
E}[Y^{2k}]={\bf E}[X^k],$ $~{\bf E}[Y^{2k-1}]$ $=0$ for
$k=1,2,\ldots$. Let $g$ be the density of $G$. Then for the
Carleman quantity of the moments of $g$ we have:
$$
{\bf C}[g]:=\sum_{k=1}^{\infty}\left({\bf E}[Y^{2k}]\right)^{-1/(2k)}
=\sum_{k=1}^{\infty}\left({\bf E}[X^{k}]\right)^{-1/(2k)}={\bf
C}[f]<\infty.
$$
This implies that for some $x_0^*>x_0$,
the logarithmic normalized integral (called also Krein quantity of $g$) over the domain $\{x:
|x|\geq x_0^*\}$ is finite:
$$
{\bf K}[g]:=\int_{|x|\geq x_0^*}\frac{-\log g(x)}{1+x^2}dx<\infty,
$$
as shown in the proof of Theorem 2 in 
\cite{l1997}. Finally,
according to Theorem 2.2 in \cite{p1998},
this is a sufficient condition for $Y$ to be M-indet on ${\mathbb
R}$ and we conclude that $X$ is M-indet on ${\mathbb R}^+$ by
mimicking the proof of Corollary 1 in \cite{s1993} (see also
\cite{p2001}, Proposition 1 and Theorem 3). \hspace{\fill} $\Box$

\vspace{0.5cm}\noindent {\bf 8. Generalized gamma distributions. Part (c) }

\vspace{0.1cm} Let us see how Theorem 4 in Section 7 works for
a random variable $\xi\sim GG(\alpha, \beta, \gamma).$
We claim that for $n>2\beta$, the power $X_n=\xi^n$ is M-indet.
To see this, recall that
$$\frac{{\bf E}[X_n^{k+1}]}{{\bf E}[X_n^{k}]}\approx
(n/\alpha\beta)^{n/\beta} (k+1)^{n/\beta}~~ \hbox{as}~
k\rightarrow \infty,$$
 where
$n/\beta>2$. Thus the moments of $X_n$ grow at a rate more than 2. Let us check that
the density $h$ of $X_n$ satisfies the
condition (2). Indeed, we have
$$
L_{h}(x):=-\frac{xh^{\prime}(x)}{h(x)}=1-\frac{\gamma}{n}+\frac{\alpha\beta}{n}
x^{\beta/n} \nearrow \infty~~~\hbox{ultimately as}~ x\rightarrow
\infty.
$$
Therefore, for $n>2\beta$, $X_n$ is M-indet by Theorem 4.

\vspace{0.2cm} \noindent {\bf Remark 3} \ To use Theorem 4 is
another way to prove some known facts, for example, that the
log-normal distribution and  the cube of the exponential
distribution are M-indet. Indeed, for $X\sim LogN(0,1),$ we have
the moment recurrence
$$
m_{k+1}=e^{k+1/2}m_k,~k=1,2,\ldots,
$$
and for $X=\xi^3$, where
$\xi\sim Exp (1)$, we have
$$m_{k+1}=(3k+1)(3k+2)(3k+3)m_k, ~k=1,2,\ldots.$$
It is easily seen that in both cases the growth rates of the
moments are quite fast. For the cube of $Exp(1)$ we have
$m_{k+1}/m_k\geq C(k+1)^3,~k=1,2,\ldots$, for some constant $C>0$,
so the rate is more than 2. For $LogN$ the rate is exponential,
hence much larger than 2. It remains to check that condition (2)
is satisfied for the density of $\xi^3$ and the density of $LogN.$
Details are omitted.

We can make one step more by considering the logarithmic
skew-normal distributions with density
$f_{\lambda}(x)=(2/x)\varphi(\ln x)\Phi(\lambda \ln x),$ $~x>0$,
where $\lambda$ is a real number. (When $\lambda=0$, $f_{\lambda}$
reduces to the standard log-normal density.) Then we have the
moment relationship
$$m_{k+1}\approx e^{(k+1/2)\rho}m_k,~~\hbox{as}~~k\rightarrow \infty,$$
where $\rho\in (0,1]$ is a constant (see, e.g.,
\cite{ls2009}, Proposition 3). Thus the moments grow very fast,
exponentially, and it remains to check that the density function
$f_{\lambda}$ satisfies the condition (2):
$$L_{f_{\lambda}}(x):=-\frac{xf_{\lambda}^{\prime}(x)}{f_{\lambda}(x)}\nearrow
\infty~~\hbox{ultimately as}~~x\rightarrow \infty.$$ Therefore, by
the above Theorem 4, we conclude that all logarithmic skew-normal
distributions are M-indet. This is one of the results in
\cite{ls2009} where a different proof is given.

\vspace{0.4cm} \noindent {\bf 9. General result on the M-indet
property of the product $Y_n=\xi_1\xi_2 \cdots \xi_n$}

\vspace{0.1cm}
 In the next theorem we describe conditions on the
distribution of $\xi$ under which the product $Y_n=\xi_1\xi_2
\cdots \xi_n$ is M-indet.

\vspace{0.2cm}\noindent {\bf Theorem 5} \ {\it Let $\xi \sim F,$
where $F$ is absolutely continuous
 with density $f>0$ on ${\mathbb R}^+$. Assume further that:}\\
(i) {\it $f(x)$ is decreasing in $x\geq 0,$ and} \\
(ii) {\it there exist two constants $x_0\geq 1$ and $A>0$ such
that
\begin{eqnarray}
f(x)/\overline{F}(x)\geq A/x~~\hbox{for}~~x\geq x_0,
\end{eqnarray} and some constants $B>0,~ \alpha>0,$
$\beta>0$ and a real $\gamma$ such that
\begin{eqnarray}\overline{F}(x)\geq Bx^{\gamma}e^{-\alpha
x^{\beta}}~~ \hbox{for}~~ x\geq x_0. \end{eqnarray} Then,  for
$n>2\beta$, the product $Y_n$ has a finite Krein quantity and is M-indet.}

\vspace{0.2cm}\noindent {\bf Corollary 2} \ {\it Let $\xi \sim F$
satisfy the conditions in Theorem 5 with $\beta<\frac12$. Then $F$
itself is M-indet.}

\vspace{0.2cm}\noindent {\bf Lemma 4} \ {\it Under the condition}
(3), {\it we have
$$
\int_x^{\infty}\frac{f(u)}{u}du\geq
\frac{A}{1+A}\frac{\overline{F}(x)}{x} \ \mbox{ and } \
\overline{F}(x)\leq \frac{C}{x^A}, ~  x>x_0,~\hbox{for some
constant}~ C>0.
$$ }
\vspace{0.1cm}\noindent {\it Proof} \ Note that for $x>x_0$,
\begin{eqnarray*}
\int_x^{\infty}\frac{f(u)}{u}du
=-\int_x^{\infty}\frac{1}{u}d\overline{F}(u)=\frac{\overline{F}(x)}{x}
-\int_x^{\infty}\frac{\overline{F}(u)}{u^2}du\geq\frac{\overline{F}(x)}{x}
-\frac{1}{A}\int_x^{\infty}\frac{f(u)}{u}du.
\end{eqnarray*}
The last inequality is due to (3).  Hence
$$\left(1+\frac{1}{A}\right)\int_x^{\infty}\frac{f(u)}{u}du\geq
\frac{\overline{F}(x)}{x}.$$ On the other hand, for $x>x_0$,
\begin{eqnarray*}\log\overline{F}(x)&=&{-\int_0^xf(t)/\overline{F}(t)dt}
={-\int_0^{x_0}f(t)/\overline{F}(t)dt-\int_{x_0}^xf(t)/\overline{F}(t)dt}\\
&\equiv&C_0-\int_{x_0}^xf(t)/\overline{F}(t)dt\leq
C_0-\int_{x_0}^xA/tdt=  C_0+A\log x_0 -A\log x.
\end{eqnarray*}
Therefore, $\overline{F}(x)\leq C/x^A,~x>x_0,$ where
$C=x_0^Ae^{C_0}.$\hspace{\fill} $\Box$

\vspace{0.2cm} \noindent{\bf Remark 4} \ After deriving in Lemma 4
a lower bound for $\int_x^{\infty}(f(u)/u)du$ we have the
following upper bound  for arbitrary density $f$ on ${\mathbb
R}^+$:
$$\int_x^{\infty}\frac{f(u)}{u}du\leq
\frac{1}{x}\int_x^{\infty}{f(u)}du= \frac{\overline{F}(x)}{x},
~x>0.$$

\vspace{0.2cm} \noindent {\it Proof of Theorem 5} \  The density
$g_n$ of $Y_n$ is expressed as follows:
$$g_n(x)=\int_0^{\infty}\!\!\int_0^{\infty}\!\!\cdots\!\!
\int_0^{\infty}\frac{f(u_1)}{u_1}\frac{f(u_2)}{u_2}
\cdots\frac{f(u_{n-1})}{u_{n-1}}f\left(\frac{x}{u_1u_2\cdots
u_{n-1}}\right)du_1du_2\cdots du_{n-1}$$ for $x>0$.  Hence
$g_n(x)>0$ and decreases in $x \in(0,\infty)$. For any $a>0$, we
have
\begin{eqnarray}g_n(x)\!\!&\geq&\!\!\int_a^{\infty}\!\!
\int_a^{\infty}\!\!\cdots\!\!\int_a^{\infty}
\frac{f(u_1)}{u_1}\frac{f(u_2)}{u_2}
\cdots\frac{f(u_{n-1})}{u_{n-1}}f\left(\frac{x}{u_1u_2\cdots
u_{n-1}}\right)du_1du_2\cdots du_{n-1}\nonumber\\
&\geq&\int_a^{\infty}\!\!\int_a^{\infty}\!\!\cdots\!\!
\int_a^{\infty}\frac{f(u_1)}{u_1}\frac{f(u_2)}{u_2}
\cdots\frac{f(u_{n-1})}{u_{n-1}}f\left(\frac{x}{a^{n-1}}\right)du_1du_2\cdots
du_{n-1}\nonumber\\
&=&f\left(\frac{x}{a^{n-1}}\right)\left(\int_a^{\infty}\frac{f(u)}{u}du\right)^{n-1},
~~x>0.
\end{eqnarray}
The above second inequality follows from the monotone property of
$f$. Taking $a=x^{1/n}>x_0$, we have, by (3)--(5) and Lemma 4, that
\begin{eqnarray*}g_n(x)&\geq&
f\left(x^{1/n}\right)\left(\int_{x^{1/n}}^{\infty}
\frac{f(u)}{u}du\right)^{n-1}\geq
f\left(x^{1/n}\right)\left(\frac{A}{1+A}\frac{\overline{F}
(x^{1/n})}{x^{1/n}}\right)^{n-1}\\
&\geq&\left(\frac{A}{1+A}\right)^{n-1}x^{-(1-1/n)}
\frac{f\left(x^{1/n}\right)}{\overline{F}(x^{1/n})}
\left(\overline{F}(x^{1/n})\right)^n\\
&\geq &C_nx^{\gamma-1}e^{-n\alpha x^{\beta/n}},
\end{eqnarray*}where
$C_n=\left(\frac{A}{1+A}\right)^{n-1}AB^n$. Therefore, the Krein
quantity for $g_n$ is as follows:
\begin{eqnarray*}
{\bf K}[g_n]&=&\int_0^{\infty}\frac{-\log
g_n(x^2)}{1+x^2}dx=\int_0^{x_0^n}\frac{-\log
g_n(x^2)}{1+x^2}dx+\int_{x_0^n}^{\infty}\frac{-\log
g_n(x^2)}{1+x^2}dx\\
&\leq &\left(-\log
g_n(x_0^{2n})\right)\int_0^{x_0^n}\frac{1}{1+x^2}dx+\int_{x_0^n}^{\infty}\frac{-\log
g_n(x^2)}{1+x^2}dx<\infty~~\hbox{if}~ n>2\beta.\end{eqnarray*}
This in turn implies that $Y_n$ is M-indet for $n>2\beta$ (see,
e.g., \cite{l1997}, Theorem 3). \hspace{\fill}$\Box$

\vspace{0.5cm}\noindent {\bf 10. Generalized gamma distributions. Part (d) }

\vspace{0.1cm} Let us see how the general result from Section 9
can be used to establish the moment indeterminacy of products of
independent copies of a random variable $\xi\sim GG(\alpha, \beta,
1).$ Here $\gamma =1$ and the density is $f(x)=ce^{-\alpha
x^{\beta}}, \ x \geq 0.$

We claim that for $n>2\beta$,
the product $Y_n=\xi_1\xi_2\cdots\xi_n$ is
M-indet. To see this, note that $f(x)/\overline{F}(x)\approx
\alpha\beta x^{\beta-1}$ and $\overline{F}(x)\approx
[c/(\alpha\beta)]x^{1-\beta}e^{-\alpha x^{\beta}}$ as $x\rightarrow
\infty$. Then the density  $f$ satisfies the conditions (i)
and (ii) in Theorem 5 and hence $Y_n$ is M-indet if $n>2\beta$.

For example, if $\xi \sim Exp(1),$ then, as mentioned before, the product
$Y_n=\xi_1\xi_2 \cdots \xi_n$ is M-indet for $n\geq 3.$

If $\xi$ has the half-normal distribution, its density is
$f(x)=\sqrt{{2}/{\pi}}e^{-x^2/2},$ $~x\geq 0$, then
$Y_n=\xi_1\xi_2 \cdots \xi_n$ is M-indet for $n\geq 5$ (recall
from Section 5 that $Y_n$ is M-det for $n\le 4$). By words: The
product of two, three or four half-normal random variables is
M-det, while the product of five or more such variables is
M-indet.

\vspace{0.2cm} In summary, we have the following result about
$GG(\alpha, \beta, \gamma)$ with $\gamma=1$.

\vspace{0.2cm} \noindent{\bf Lemma 5} \ {\it Let $n\geq
2,~X_n=\xi^n$ and $Y_n=\xi_1 \cdots \xi_n$, where $\xi_1,
\ldots,\xi_n$ are independent copies of $\xi\sim GG(\alpha, \beta,
1)$. Then the power  $X_n$ is M-det iff the product $Y_n$ is M-det
and this is true iff $n\leq 2\beta$.}

\vspace{0.1cm} We now consider the general case $\gamma>0$.

\vspace{0.2cm} \noindent {\bf Theorem 6} \ {\it Let $n\geq
2,~X_n=\xi^n$ and $Y_n=\xi_1\cdots\xi_n$, where
$\xi_1,\ldots,\xi_n$ are independent copies of $\xi\sim GG(\alpha,
\beta, \gamma)$. Then $X_n$ is M-det iff $Y_n$ is M-det and this
is true iff $n\leq 2\beta$. In other words, both $X_n$ and $Y_n$
have the same moment determinacy property.}

\vspace{0.1cm}\noindent {\it Proof} \ Define $\eta=\xi^{\gamma}$,
$\eta_i=\xi_i^{\gamma}$, $i=1,2,\ldots,n$,
$X_n^*=\eta^n=(\xi^n)^{\gamma}=X_n^{\gamma}$ and
$Y_n^*=\eta_1\eta_2\cdots\eta_n=(\xi_1\xi_2\cdots\xi_n)^{\gamma}=Y_n^{\gamma}$.
Since $\eta\sim GG(\alpha, \beta/\gamma, 1)$, we have, by Lemma 5,
$X_n^*$ is M-det iff $Y_n^*$ is M-det iff $n\leq 2\beta/\gamma$.
Next, note that for each $x>0$, we have ${\bf P}[X_n^*>x]={\bf
P}[X_n>x^{1/\gamma}]$ and ${\bf P}[Y_n^*>x]={\bf
P}[Y_n>x^{1/\gamma}]$. This implies that any distributional
property shared by $X_n^*$ and $Y_n^*$ can be transferred to a
similar property shared by $X_n$ and $Y_n$, and vice versa.
Therefore, $X_n$ is M-det iff $Y_n$ is M-det iff $n\leq 2\beta$,
because $X_n$ is M-det iff $n\leq 2\beta$ (see, e.g.,
\cite{pk1992}). \hspace{\fill}$\Box$

\vspace{0.5cm} \noindent {\bf 11. Half-logistic distribution}

\vspace{0.2cm} Some of the above results and  illustrations
involve the generalized gamma distribution $GG(\alpha, \beta,
\gamma)$. It is useful to have a moment determinacy
characterization for powers and products based on non-GG
distributions. Here is an example based on the half-logistic
distribution, which clearly is not in the class $GG$.

\vspace{0.2cm} \noindent
{\bf Statement} \ {\it We say that the random variable $\xi$ has
the half-logistic distribution if its density is
$$
f(x)=\frac{2e^{-x}}{(1+e^{-x})^2},~~x\geq 0.
$$
The power $X_n=\xi^n$ and the product $Y_n=\xi_1\xi_2 \cdots
\xi_n$ are defined as above. Then $X_n$ is M-det iff $Y_n$ is
M-det and this is true iff $n\leq 2$. This means that for each
$n$, the two random variables $X_n$ and  $Y_n$ share the same
moment determinacy property.}

\vspace{0.1cm}\noindent {\it Proof} \
(i) The claim that $X_n$ is M-det iff $n\leq 2$ follows from results in 
\cite{lh1997}. Actually, in \cite{lh1997} it is proved that for any real $s>0$, the power
$\xi^s$ is M-det iff $s\leq 2$. Let us give here an alternative proof.
The density $h_s$ of $\xi^s$ is
$$h_s(z)=\frac{2}{s}z^{1/s-1}\frac{e^{-z^{1/s}}}{(1+e^{-z^{1/s}})^{2}},~z\geq 0.$$
Using the inequality: $\frac{1}{4}\leq (1+e^{-x})^{-2}\leq 1$ for $x\geq
0$, we find two-sided bounds for the moments of $\xi^s$:
$$\frac{1}{2}\Gamma(ks+1)\leq {\bf E}[(\xi^s)^k]\leq
\int_0^{\infty}\frac{2}{s}z^{k+1/s-1}{e^{-z^{1/s}}}dz=
2\Gamma(ks+1).$$ Therefore the growth rate of the moments of $\xi^s$ is
$$
\frac{{\bf E}[(\xi^s)^{k+1}]}{{\bf E}[(\xi^s)^k]}\leq
4\cdot\frac{\Gamma((k+1)s+1)}{\Gamma(ks+1)}\approx 4s^s
(k+1)^s~~\hbox{as}~k\rightarrow\infty.
$$
By Theorem 1, this implies that $\xi^s$ is M-det for $s\leq 2$.
On the other hand, we have
$$
\frac{{\bf E}[(\xi^s)^{k+1}]}{{\bf E}[(\xi^s)^k]}\geq
\frac{1}{4}\cdot\frac{\Gamma((k+1)s+1)}{\Gamma(ks+1)}\approx
\frac{1}{4}s^s (k+1)^s~~\hbox{as}~k\rightarrow\infty.
$$
The moment condition in Theorem 4 is satisfied if $s>2$. It remains now to
check the validity of condition (2) for the density $h_s$. We have
\begin{eqnarray*}&&L_{h_s}(z):=-\frac{zh_s^{\prime}(z)}{h_s(z)}\\
&=&1-\frac{1}{s}+\frac{1}{s}z^{1/s}-
\frac{2}{s}\,z^{1/s}\frac{e^{-z^{1/s}}}{1+e^{-z^{1/s}}}~~\nearrow
\infty~\hbox{ultimately~~as}~~z\rightarrow \infty.\end{eqnarray*}
Hence, if $s>2$, ~$\xi^s$ is M-indet.\\
(ii) It remains to prove that $Y_n$ is M-det iff $n\leq 2$. \\
(Sufficiency) As in part (i),  we have
$$\frac{1}{2}\Gamma(k+1)\leq {\bf E}[\xi^k]= 2\Gamma(k+1).$$ Therefore,
${\bf E}[\xi^{k+1}]/{\bf
E}[\xi^{k}]={\cal O}(k+1)$ as $k\rightarrow\infty$. By
Theorem 2, we conclude that $Y_n$ is M-det if $n\leq 2$. \\
(Necessity) Note that $\overline{F}(x)={\bf P}[\xi>x]=
2e^{-x}/(1+e^{-x})\geq e^{-x},~~x\geq 0,$ and
$f(x)/\overline{F}(x)=1/(1+e^{-x})\geq 1/2,~x\geq 0$. Therefore,
taking $\beta=1$ in Theorem 5, we conclude that $Y_n$ is M-indet
if $n> 2$. Let us express this conclusion by words: The product of
three or more half-logistic random variables is
M-indet.\hspace{\fill} $\Box$

\vspace{0.5cm}\noindent {\bf Acknowledgments}

\vspace{0.1cm} We would like to thank the Editors and the Referees
for their valuable suggestions which helped us shorten the
original manuscript and improve the presentation. We thank also
Prof. M.C. Jones (Open University, UK) for showing interest in our
work and making useful comments.

This paper was basically completed during the visit of JS to
Academia Sinica, Taipei, December 2012 -- January 2013. The
support and the hospitality provided by the Institute of
Statistical Science are greatly appreciated.

The work of GDL was
partly supported by the National Science Council of ROC (Taiwan)
under Grant  NSC 102-2118-M-001-008-MY2, and that of JS by
the Emeritus Fellowship provided by the Leverhulme Trust (UK).

\vspace{0.2cm}
\bibliographystyle{plain}


\end{document}